\documentclass[11pt]{article}
\usepackage{enumerate}
\usepackage{lscape}
\usepackage{appendix}
\usepackage{xcolor}
\usepackage{changepage}
\usepackage{amssymb,latexsym}
\usepackage{graphicx}
\usepackage{fancyhdr}
\usepackage{mathdots}
\usepackage{enumitem}
\usepackage{cite}
\usepackage[english]{babel}
\usepackage{bbm}
\usepackage{amssymb}
\usepackage{tabularx}
\usepackage{epsfig}
\usepackage{amsfonts}
\usepackage{amscd}
\usepackage{amsmath,amsthm}
\usepackage{lmodern}
\usepackage{mathrsfs}
\usepackage{url}
 \usepackage{arydshln}
\usepackage[linktocpage, colorlinks=true,citecolor=blue,linkcolor=blue,urlcolor=blue]{hyperref}

\numberwithin{equation}{section}
\theoremstyle{plain}

\newtheorem{theorem}{Theorem}[section]

\newtheorem{proposition}[theorem]{Proposition}

\theoremstyle{definition}
\newtheorem{definition}[theorem]{Definition}
\newtheorem{example}[theorem]{Example}

\newtheorem{identity}{Identity}
\newtheorem{proplem}{Proposition-Lemma}[section]

\theoremstyle{remark}
\newtheorem{remark}[theorem]{Remark}

\renewcommand\appendix{\par
  \setcounter{section}{0}
  \setcounter{subsection}{0}
  \setcounter{figure}{0}
  \setcounter{table}{0}
  \renewcommand\thesection{Appendix \Alph{section}}
  \renewcommand\thefigure{\Alph{section}\arabic{figure}}
  \renewcommand\thetable{\Alph{section}\arabic{table}}
}
\newcommand{\be}{\begin{equation}}
\newcommand{\ee}{\end{equation}}
\newcommand{\bea}{\begin{eqnarray}}
\newcommand{\eea}{\end{eqnarray}}
\newcommand{\bid}{\begin{identity}}
\newcommand{\eid}{\end{identity}}
\newcommand{\bin}{\begin{inequality}}
\newcommand{\ein}{\end{inequality}}
\newcommand{\bp}{\begin{proof}}
\newcommand{\ep}{\end{proof}}

\newcommand{\Z}{{\mathbb Z}}

\newcommand{\enom}[3]{{#1 \choose #2}_{\hskip -3pt #3}}

\def\T {{\mathsf{T}_{\hskip-1pt m}}}

\title{  \bf Some identities involving \\polynomial coefficients}

\author{Nour-Eddine Fahssi\\[0.8ex]
\small \it Department of Mathematics, FSTM, \\
\small \it Hassan Second University of Casablanca,\\
\small \it BP 146, Mohammedia,  Morocco.\\
\small \href{mailto:n.fahssi@live.fr}{\tt n.fahssi@live.fr}}

\date{}
\begin{document}
\maketitle


\begin{abstract}
By polynomial (or extended binomial) coefficients, we mean the coefficients in the expansion of integral powers, positive and negative, of the polynomial $1+t +\cdots +t^{m}$; $m\geq 1$ being a fixed integer.  We will establish several identities and summation formul\ae\ parallel to those of the usual binomial coefficients.\\
[0.8ex]
\small \noindent 2010 AMS Classification Numbers: 11B65, 05A10, 05A19.
\end{abstract}
\section{Introduction and preliminaries}
\begin{definition} Let $m \geq 1$ be an integer and let $p_m(t) = 1+t +\cdots +t^{m}$. The coefficients defined by \footnote{We make use of the conventional notation for coefficients of entire series : $\left[t^n\right]\sum_i a_i t^i~:= a_n$} \be \label{enom}  \enom{n}{k}{m} \stackrel{\mbox{\tiny def}}=  \left\{
                                                     \begin{array}{ll}
                                                       [t^k]\left(p_m(t)\right)^n, &  \hbox{if} \quad k \in \Z_{\geq }, \; n \in \Z \\
                                                       0, & \hbox{if} \quad k \not\in \Z_{\geq } \; \; \hbox{or} \; \; k > mn \;\, \hbox{whenever} \; n \in \Z_{> } ,                                                     \end{array}
                                                   \right.\ee are called \emph{polynomial coefficients}~\cite{Comt} or \emph{extended binomial coefficients}. When the rows  are indexed by $n$ and the columns by $k$, the polynomial coefficients form an array called \emph{polynomial array} of degree $m$ and denoted by $\T$. The sub-array corresponding to positive $n$ is termed \emph{extended Pascal triangle} of degree $m$, or \emph{Pascal-De Moivre triangle}, and denoted by $\mathsf{T}_m^+$. For instance, the array $\mathsf{T}_3$ begins as shown in Table~\ref{fig}.\end{definition}
The study of the coefficients of $\mathsf{T}_m^+$ is an old problem. Its history dates back at least to De Moivre~\cite{moivre,balasu} and Euler~\cite{Euler}.  In modern times, they have been studied in many papers, including several in The Fibonacci Quarterly. See~\cite{fahssi} for a systematic treatise on the subject and an extensive bibliography.

As well as of being interest for their own right, the polynomial coefficients have applications in number theory and combinatorics~\cite{fahssi,Freund}. For that reason, the search for identities involving them is important. This paper presents several identities, summations and generating functions that generalize those fulfilled by the binomial coefficients. The author believes that most of the results in this paper are new.

In section~\ref{s2} we will prove extensions of some of most well-known binomial coefficient identities (Table~\ref{tab1}). In section~\ref{s3} we provide two types of generating functions (GF). To go further in the extension, we prove in section~\ref{s4} a dozen of identities, and give five non-trivial summation formulas involving trinomial ($m=2$) coefficients. In section \ref{s5}, we conclude with some remarks and suggestions for further research.

It is well-known that (see~\cite{Freund} and~\cite[p.104]{riordan1})
 \bea \label{4term} \enom{n}{k}{m}&=&\enom{n}{k-1}{m}+ \enom{n-1}{k}{m}-\enom{n-1}{k-m-1}{m},\\ \label{Tbis} \enom{n}{k}{m}&=&\sum_{i=\lceil k/m \rceil}^k  \enom{i}{k-i}{m-1} {n \choose i}. \eea where $\lfloor \cdot \rfloor$ and $\lceil \cdot \rceil$ denote, respectively, the floor and ceiling functions. By using the upper negation property of binomial coefficients : \mbox{${-n \choose k}= (-1)^k \binom{n+k-1}{k}$}, it is easy to check that the four-term recurrence~\eqref{4term} and the equation~\eqref{Tbis} hold also true for negative $n$.

The sequence $\{\enom{-1}{k}{\,m}\}_{_k}$ will be involved in many identities. For notational convenience, we denote it by $\chi_m(k)$. It is easy to show that  $\chi_0(k)=[k=0]$ \footnote{Throughout the paper we use the Iverson bracket notation to indicate that $[P]=1$ if the proposition $P$ is true and 0 otherwise.}, $\chi_1(k)=(-1)^k$, and generally
\be \chi_m(k)= \left\{\begin{array}{ll}
1 & \hbox{if} \quad k \equiv 0 \pmod{m+1} \\
-1 & \hbox{if} \quad k \equiv 1 \pmod{m+1}  \\
0 & \hbox{otherwise.}
\end{array}\right.
\ee
\begin{table}
\centering
\[\begin{CD} @ >{k}>>\\
\end{CD}\]
\[\begin{CD}
@V{n}VV \\
\end{CD} \begin{array}{c||ccccccccccc}
&0&1&2&3&4&5&6&7&8&9&\ldots\\
\hline \hline
\vdots&\vdots &\vdots &\vdots &\vdots &\vdots &\vdots &\vdots  & \vdots &\vdots &\vdots &\iddots \\
-3& 1 & -3 & 3 & -1 & 3 & -9 & 9  & -3 & 6 & -18 & \ldots \\
-2& 1 & -2 & 1 & 0 & 2 & -4 & 2 & 0 &3 &-6 & \ldots\\
-1& 1 & -1 & 0 & 0 & 1 & -1 & 0 &0&1&-1&\ldots \\ \hdashline
0& 1 &  &   &  &  &  &  &&&&\\
1& 1 & 1 & 1 & 1  &  &  &&&&&\\
2& 1 & 2 & 3 & 4 & 3 & 2 & 1 &&&& \\
3& 1 & 3 & 6 & 10 & 12 & 12 & 10 & 6&3&1\\
\vdots&\vdots &\vdots &\vdots &\vdots &\vdots &\vdots &\vdots & \vdots&\vdots&\vdots&\ddots
\end{array}
\]
\caption{The beginning of the array $\mathsf{T}_{\hskip-1pt 3}$}\label{fig}
\end{table}
\section{Extension of most known binomial coefficient identities} \label{s2}
\begin{proplem}The extended binomial coefficient identities in Table~\ref{tab1} hold true.\end{proplem}
\bp (sketch) The factorial expansion (i) for positive $n$ and finite $m$ is known~\cite{Bonda}; it is established by an application of the usual multinomial theorem; the sum runs over all $m$-tuples such that $\sum_{i \leq m}i n_i =k$.  {\small \begin{landscape}
\begin{table}[htbp]
\thispagestyle{empty}
{{\renewcommand{\arraystretch}{0.7}
\begin{tabular}{llllll}
\hline \hline
&\textbf{Identity}  & \multicolumn{1}{c}{\textbf{Binomial\cite{gould}}} &&\multicolumn{1}{c}{\textbf{Polynomial} }& \\
\hline &&&&  &\\
(i) & Factorial expansion &  $\displaystyle{n \choose k}=\frac{n!}{k!(n-k)!}$ && $\displaystyle \enom{n}{k}{m}= \sum_{n_i \leq n}  \frac{n!}{(n-n_1- \ldots - n_m)!\, n_1! \cdots n_m!} \Big[\sum_{i \leq m} i n_i = k\Big] $ & integer $n \geq 0$\\
&&&&&\\
(ii)& Symmetry &$\displaystyle\binom{n}{k}=\binom{n}{n-k}$ && $\displaystyle \enom{n}{k}{m}=\enom{n}{m n-k}{m}$ & integer $n \geq 0$ \\
& &&&& \\
(iii) & Absorption/Extraction  &$\displaystyle \binom{n}{k}=\frac{n}{k}\binom{n-1}{k-1}$  &&$ \displaystyle \enom{n}{k}{m} = \frac{n}{k} \sum_{i=1}^m i \enom{n-1}{k-i}{m}$& integer $n$; $k \neq 0$\\
& && &&\\
(iv)&  ${\mbox{Vandermonde} \atop \mbox{convolution}}$
   &$ \displaystyle \sum_{i+j=k}\binom{r}{i}
\binom{s}{j}= \binom{r+s}{k}$ && $\displaystyle
\sum_{i+j=k}\enom{r}{i}{m}
\enom{s}{j}{m}= \enom{r+s}{k}{m} $ & integers $r$ and $s$\\
& && &&\\
(v)& Addition/Induction  &$\displaystyle \binom{n}{k}=\binom{n-1}{k}+\binom{n-1}{k-1}$  && $\displaystyle \enom{n}{k}{m} = \sum_{i=0}^{m} \enom{n-1}{k-i}{m}$ & integer $n$\\
& & &&&\\
(vi)& Binomial theorem  &$\displaystyle \sum_{k \geq 0} \binom{n}{k} x^k y^{n-k}=(x+y)^n$  && $\displaystyle \sum_{k \geq 0} \enom{n}{k}{m} x^k y^{m n-k}=\left(\sum_{i=0}^m  x^i y^{m-i}\right)^n $ & $\displaystyle {\left| p_m(x/y)-1 \right| < 1 \atop \mbox{if $n<0$ }}$\\
& &&& &\\
(vii) & Upper summation &$\displaystyle \sum_{0 \leq l \leq n}{l \choose k}={n+1 \choose k+1}$ && $\displaystyle \sum_{0 \leq l \leq n} \enom{l}{k}{m}=\sum_{i=0}^k \chi_{m-1}(i) \enom{n+1}{k-i+1}{m}$& integer $n \geq 0$\\
& &&&&\\
(viii)& Parallel summation &$\displaystyle \sum_{0 \leq k \leq n}{r +k \choose k}={r+n+1 \choose n}$ && $\displaystyle  \sum_{0 \leq k \leq n} \enom{r+k}{mk}{m} = \sum_{i=0}^{mr}\chi_{m-1}(i)\enom{r+n+1}{m r-i+1}{m}$& integer $r \geq 0$\\
&&&&&\\
(ix)& Horizontal recurrence & $\displaystyle k {n \choose k}=\left(n+1-k\right) {n \choose k-1}$ && $\displaystyle  k \enom{n}{k}{m}=\sum_{i=1}^{m}\left((n+1)i-k\right) \enom{n}{k-i}{m}$&  integer $n$\\
&&&&&\\
\hline \hline
\end{tabular}}}
\medskip
\caption {Extensions of the most important binomial coefficient identities.}
\label{tab1}
 \end{table}
\end{landscape} } \noindent The polynomial symmetry (ii) is a consequence of
the fact that the polynomial $p_m$ is  self-reciprocal : $p_m(t)=t^m p_m (t^{-1})$. The absorption/extraction identity (iii) follows from taking the derivatives of both sides of $p_m(t)^n=\sum_{k \geq 0} \enom{n}{k}{\, m}t^k$, and equating the coefficients of $t^k$. The Vandermonde convolution (iv) is obtained by equating the coefficients of $t^k$ in the sides of $p_m^{r+s}(t)=p_m^r(t) p_m^s(t)$. The addition/induction (v) is a particular case of Vandermonde convolution with $r=1$ and $s=n-1$. The generalized binomial theorem (vi) is an obvious consequence of the definition~\eqref{enom}. We prove the upper summation (vii) as follows : \begin{multline*} \sum_{0 \leq l \leq n} \enom{l}{k}{m} = [t^k] \sum_{0 \leq l \leq n} (p_m(t))^l = [t^k] \frac{p_m(t)^{n+1}-1}{p_m(t)-1}=[t^{k}] \frac{1}{p_{m-1}(t)} \frac{p_m(t)^{n+1}-1}{t} =\\
 \sum_{i=0}^k [t^i]\frac{1}{p_{m-1}(t)} \, [t^{k-i}]\frac{p_m(t)^{n+1}-1}{t}
= \sum_{i=0}^k \chi_{m-1}(i) \enom{n+1}{k-i+1}{m}. \end{multline*} \noindent To prove the parallel summation (viii), we write
$\sum_{0 \leq k \leq n} \enom{r+k}{mk}{\, m} = \sum_{l=r}^{r+n} \enom{l}{mr}{\, m} = \sum_{l=0}^{r+n} \enom{l}{mr}{\, m}$, and then, apply  Identity (vii) and the symmetry (ii). To prove the horizontal recurrence (ix), we combine Identity (iii) with Identity (v).  \ep
\section{Generating functions} \label{s3}
Let \[ F_k^+(x)=\sum_{n=0}^\infty \enom{n}{k}{m}x^n, \quad \hbox{and} \quad F_k^-(x)=\sum_{n=1}^\infty \enom{-n}{k}{m}x^n \] be, respectively, the series generating the $k$-th column of the triangle $\mathsf{T}_m^+$ and the $k$-th column of the upper negated sub-array $\mathsf{T}_m^-$.
\begin{theorem}[Vertical generating functions] \label{cgf}
The functions $F_k^+$ and $F_k^-$ take the forms \be \label{GF}
F_k^+(x)=\left(\frac{1} {1-x}\right)^{k+1} P_k^{(m)}(x),  \quad \hbox{and} \quad F_k^-(x)=(-1)^k \left(\frac{x} {1-x}\right)^{k+1} P_k^{(m)}(x^{-1})\, ,\ee where $ P_k^{(m)}$ are polynomials generated by
 \be \label{GFb} \sum_{k=0}^\infty P_k^{(m)}(x) y^k = \frac{1+(x-1) y}{1-y+x (1-x)^m y^{m+1}}\, , \ee
 or, alternatively, derived from the recurrence relation \be \label{recur} P_k^{(m)}(x)=P_{k-1}^{(m)}(x)-x (1-x)^m P_{k-m-1}^{(m)}(x),\ee with  $P_0^{(m)}(x)=1$, $P_1^{(m)}(x)=x$ and $P_k^{(m)}(x)=0$ for $k <0$. In particular, $P_k^{(m)}(x)=x$ for all $1 \leq k \leq m$.
\end{theorem}
\bp
Using the recurrence~\eqref{Tbis}, we write \[F_k^+(x)= \sum_{i=\lceil k/m \rceil}^k \enom{i}{k-i}{m-1}
\sum_{n=0}^\infty {n \choose i} x^n.\] So, from the GF of binomial coefficients~:~$\sum_{n=0}^\infty {n \choose
k}u^n = \frac{u^k}{(1-u)^{k+1}}$, the function $F_k^+(x)$ can be
displayed in the form~\eqref{GF} with \be \label{expli} P_k^{(m)}(x)= \sum_{i=\lceil k/m \rceil}^k \enom{i}{k-i}{m-1} x^i(1-x)^{k-i}=\sum_{i=\lceil k/m \rceil}^k \alpha_i (k) x^i,\ee where $\alpha_i (k)$ are computable coefficients. Moreover, \[\sum_{k=0}^\infty P_k^{(m)}(x) y^k = \sum_{k=0}^\infty (1-x)^{k+1}F_k^+(x) y^k= \frac{1-x}{1-x \, p_m((1-x)y)}= \frac{1+(x-1) y}{1-y+x (1-x)^m y^{m+1}},\] where we have used the fact that $\sum_{n,k} \enom{n}{k}{m} x^n y^k = \big(1-x\,p_m(y)\big)^{-1}$. As for the recurrence relation~\eqref{recur},  we use the GF of Eq.~\eqref{4term} with respect to $n$ to get \begin{eqnarray*}
F_k^+(x)&=&F_{k-1}^+(x) +\chi_m(k)-\chi_m(k-m-1)+x F_k^+(x) -x F_{k-m-1}^+(x)\\ &=&F_{k-1}^+(x) +x F_k^+(x) -x F_{k-m-1}^+(x),
\end{eqnarray*} because $\chi_m(k)=\chi_m (k-m-1)$ for all $k$. So, $(1-x)F_k^+(x)= F_{k-1}^+(x) -x F_{k-m-1}^+(x)$, with $F_0^+(x)=(1-x)^{-1}$, and $F_1^+(x)=x (1-x)^{-2}$. Multiplying both sides of the last recurrence by $(1-x)^k$, we get the desired result. The expression of $F_{k}^-(x)$ follows in the same way by employing the formula~:~$\sum_{n=0}^\infty {-n \choose i} u^n = (-1)^i\frac{ u}{(1-u)^{i+1}}$.
\ep
\begin{example}In the special case $m=2$, the polynomials  $P_k^{(2)}$ can be explicitly determined by partial fraction decomposition of~\eqref{GFb}. This yields
\be \label{pm2} P_k^{(2)}(x)=\frac{\left(x+\sqrt{x(4-3 x)
   }\right)^{k+1}-\left(x-\sqrt{x(4-3 x)
   }\right)^{k+1}}{2^{k+1}\sqrt{x(4-3 x)}}\, .\ee\end{example}
A different generating function is given by
\begin{proposition}[Carlitz Generating Functions] For any integers $a$ and $b$, we have the following GF \be G_m(a,b;x):=\sum_{k=0}^\infty \enom{a+b k}{k}{m} x^k= \frac{\left(p_m (y)\right)^{a+1}}{p_m(y)-b \,y \,p\, '_m(y)},\ee where the indeterminates $x$ and $y$ are related by $y=x \left(p_m(y)\right)^b$.  \end{proposition}\bp The proposition is a special case of a theorem of Carlitz~\cite{carlitz} on generating functions of the form $A(x)B(x)^n$, where $A(x)$ and $B(x)$  are formal power series satisfying $A(0)=B(0)=1$. \ep \begin{remark} Closed-form expressions for Carlitz GF can be obtained for $-3 \leq mb \leq 4$. In particular, we recover the known Euler GF for central trinomial coefficients: \[ G_2 (0,1;x)=\frac{1}{\sqrt{(1+x)(1-3x)}}.\] \end{remark}
\section{Some extended combinatorial identities}\label{s4}
Now we prove more identities involving polynomial coefficients.
\begin{identity}[Extended entry 1.5 in \cite{gould}] For $n \in \Z$,  \[\sum_{j=0}^k \chi_m(j) \enom{n}{k-j}{m} = \enom{n-1}{k}{m} \, .\] \end{identity}
\bp The identity follows from the definition of $\chi_m(j)$ and the Vandermonde convolution. \ep
\begin{identity}[Extended entry 1.23 in \cite{gould}] For $n \geq 0$,  \[\sum_{{k=0 \atop m | n+k}}^\infty \enom{(n+k)/m}{k}{m} 2^{-k/m} =2^{1+\frac{n}{m}} f_{n}^{(m)}, \] where $f_n^{(m)}$ are integers satisfying the recurrence $f_n^{(m)}=2 f_{n-1}^{(m)} -f_{n-m-1}^{(m)}$ with $f_0^{(m)}=1$ and $f_n^{(m)}=0$ whenever $n < 0$.\eid \bp We have, from equations~\eqref{GF} and~\eqref{expli}, \[ \sum_{{k=0 \atop m | n+k}}^\infty \enom{(n+k)/m}{k}{m} 2^{-k/m} = \sum_{l=\lceil n/m \rceil}^\infty \enom{l}{n}{m} 2^{-l+n/m}=2^{n/m}F_n^+(1/2)= 2^{1+\frac{n}{m}}f_{n}^{(m)},\] where $f_{n}^{(m)} = 2^n P_n^{(m)}(1/2)$. The recurrence  $f_n^{(m)}=2 f_{n-1}^{(m)} -f_{n-m-1}^{(m)}$ follows from~\eqref{recur}.\ep \begin{remark} $f_n^{(2)}= F_{n+1}$ is the $(n+1)$-th Fibonacci number. \end{remark}
\begin{identity}[Extended Identity 175 in \cite{benj}] \label{175} For $n \geq 0$, \[\sum_{k \leq n}(-1)^{n-k} \enom{n-k}{k}{m} = \chi_{m+1}(n).\] \eid \bp The proof is most easily carried out by using the Herb Wilf's Snake-Oil~\cite{wilf}. We will prove that the GF of the sum in LHS is equal to $1/p_{m+1}(x)$:  \begin{eqnarray*}  \sum_{k=0}^\infty \sum_{n=k}^\infty (-1)^{n-k} \enom{n-k}{k}{m} x^n &=&\sum_{k=0}^\infty x^k \sum_{i=0}^\infty (-1)^i \enom{i}{k}{m} x^{i} \\ &\stackrel{\mbox{\tiny \eqref{GF}} }=& \sum_{k=0}^\infty x^k \left(\frac{1}{1+x}\right)^{k+1}P_k^{(m)}(-x) \\  &\stackrel{\mbox{\tiny \eqref{GFb}} }=& \dfrac{1}{1+x} \dfrac{1+(-x-1)\dfrac{x}{1+x}}{1-\dfrac{x}{1+x}-x(1+x)^m \dfrac{x^{m+1}}{(1+x)^{m+1}}} \\ &=& \dfrac{1-x}{1-x^{m+2}} = \dfrac{1}{p_{m+1}(x)}.\end{eqnarray*} \ep
\begin{identity}[Extended entry 1.68 in \cite{gould}] \label{4}  For $n \geq 1$, \[ \arraycolsep=1pt\def\arraystretch{2.5}
\sum_{k =0}^{\lfloor mn/(m+1) \rfloor} (-1)^{n-k} \enom{n-k}{k}{m}\frac{1}{n-k} = \left\{
                                  \begin{array}{lll}
                                    \displaystyle {m+1 \over n}, & & \hbox{if} \;\; m+2 \;|\; n \\
                                    \displaystyle -{1 \over n}, & & \hbox{otherwise.}
                                  \end{array}
                                \right. \] \eid \bp This identity can be cast in the following compact form : \be \label{id5} \sum_{k =0}^{\lfloor mn/(m+1) \rfloor} (-1)^{n-k} \enom{n-k}{k}{m}\frac{n}{n-k} = (m+1)[m+2 \, |\, n]-[m+2 \not| \, n], \quad n \geq 1. \ee  First, we have by application of the identities (iii). and \ref{175} : \begin{eqnarray*} && \sum_{k =0}^{\lfloor mn/(m+1) \rfloor} (-1)^{n-k} \enom{n-k}{k}{m}\left( 1+\frac{k}{n-k}\right) \\ 
                                	& \stackrel{\mbox{\tiny Id. \ref{175}} }=& \chi_{m+1}(n)+\sum_{k =0}^{\lfloor mn/(m+1) \rfloor} (-1)^{n-k} \enom{n-k}{k}{m}\frac{k}{n-k}\\ 
                                	&\stackrel{\mbox{\tiny Id. (iii)}}=&   \chi_{m+1}(n)+ \nonumber \sum_{k =0}^{\lfloor mn/(m+1) \rfloor} (-1)^{n-k} \sum_{i=1}^m i \enom{n-k-1}{k-i}{m}\\ 
                                	&=& \chi_{m+1}(n)+\sum_{i=1}^m i \sum_{k =i}^{\lfloor mn/(m+1) \rfloor} (-1)^{n-k} \enom{n-k-1}{k-i}{m}
\\ &=& \chi_{m+1}(n)+  \sum_{i=1}^m i \sum_{l =0}^{\lfloor mn/(m+1) \rfloor-i} (-1)^{n-l-i} \enom{n-l-i-1}{l}{m} \\ &\stackrel{\mbox{\tiny Id. \ref{175}}}=&  \chi_{m+1}(n) - \sum_{i=1}^m i \chi_{m+1}(n-1-i). \end{eqnarray*}
Now, it is easy to show that both the GF of the last term (with respect to $n$) and the GF of RHS of Equation~\eqref{id5} are given by : \[-1+\dfrac{1-x^2 p\,'_{m+1}(x)}{p_{m+1}(x)}=-\frac{(m+2) x^{m+2}}{x^{m+2}-1}-\frac{x}{1-x}.\]  This completes the proof of Identity~\ref{4}.\ep
\begin{identity}[Extended Entries 1.90 and 1.96 in \cite{gould}] For $n \geq 1$, \[\begin{array}{ccc}
                                                                               \displaystyle \sum_{k=0}^{\lfloor mn/2 \rfloor} (-1)^k \enom{n}{2k}{m}  &=&\left\{
\begin{array}{ll}
                                                            1,  & \hbox{if} \;\; m \equiv 0 \pmod{4} \\
                                                            2^{n/2}\cos \, (n\pi/4),  & \hbox{if} \;\; m \equiv 1 \pmod{4}\\
                                                            \cos \, (n\pi/2),  & \hbox{if} \;\; m \equiv 2 \pmod{4} \\
                                                            0 , & \hbox{if} \;\; m \equiv 3 \pmod{4}.
                                                          \end{array}
                                                        \right. \\
                                                        & & \\
             \displaystyle \sum_{k=0}^{\lfloor (mn-1)/2 \rfloor} (-1)^k \enom{n}{2k+1}{m} &=& \left\{
                                                          \begin{array}{ll}
                                                              0,& \hbox{if} \;\; m \equiv 0 \pmod{4} \\
                                                             2^{n/2}\sin \, (n\pi/4) ,& \hbox{if} \;\; m \equiv 1 \pmod{4}\\
                                                              \sin \, (n\pi/2) ,& \hbox{if} \;\; m \equiv 2 \pmod{4} \\
                                                             0,& \hbox{if} \;\; m \equiv 3 \pmod{4}.
                                                          \end{array}
                                                        \right.
                                                                             \end{array}
 \] \eid \bp We observe that the two alternating sums are the real and the imaginary parts of $ \left(p_m(\mathbf{i})\right)^n $, ($\mathbf{i}=\sqrt{-1}$). To establish them, we make use of the following formul\ae\, which can be easily shown by induction on $n$ \bea \nonumber \Re \left(1-{\rm \mathbf{i}}^{m+1}\right)^n &=& \left\{
                                                           \begin{array}{ll}
                                                             2^{n-1}(1+(-1)^p), & \hbox{if} \;\; m=2p+1 \\
                                                             2^{n/2}\cos(n\pi/4), & \hbox{if} \;\; m=2p
                                                           \end{array}
                                                         \right.\\
\nonumber \Im\left(1-{\rm \mathbf{i}}^{m+1}\right)^n &=& \left\{
                                                           \begin{array}{ll}
                                                             0, & \hbox{if} \;\; m=2p+1 \\
                                                             2^{n/2}(-1)^{p+1} \sin(n\pi/4), & \hbox{if} \;\; m=2p
                                                           \end{array}
                                                         \right.
\eea \ep
\begin{identity}[Extended Identity 5.23, Table 169 in~\cite{graham}] \label{id4} For $r, s \geq 0$ \[  \sum_{l} \enom{r}{q+l}{m}
\enom{s}{k+l}{m} = \enom{r+s}{m r-q+k}{m}=\enom{r+s}{m s+q-k}{m}, \quad r, s \geq 0 .\]
\eid \bp This identity follows from the symmetry relation (ii) and the application of Vandermonde convolution (iv). \ep
\begin{identity}[Extended Entries 3.78 and 3.79 in~\cite{gould}] \label{squares2} For all $n \geq 0$,
\begin{eqnarray*}
\sum_{k=0}^{mn} \enom{n}{k}{m}^2 &=& \enom{2n}{m n}{m}\\
 \sum_{k=0}^{mn} k \enom{n}{k}{m}^2  &=& \frac{mn}{2} \enom{2n}{mn}{m} \\
\sum_{k=0}^{mn} k^2 \enom{n}{k}{m}^2 &=& \frac{n^2}{2n-1} \sum_{i=1}^m i(m(n-1)+i) \enom{2n-1}{mn-i}{m}.
\end{eqnarray*} \eid \bp Putting $r=s=n$ and $q=k=0$ in Identity~\ref{id4} we get the first equality. The LHS of the second equation can be written as { \begin{eqnarray*}   n \sum_{k =0}^{mn} {k \over n} \enom{n}{k}{m} \enom{n}{k}{m} &\stackrel{\mbox{\tiny Id. (iii)} }=& n \sum_{i=1}^m i \sum_{k =i}^{mn} \enom{n-1}{k-i}{m}\enom{n}{k}{m} \\ &=& n \sum_{i=1}^m i \sum_{l =0}^{mn-i} \enom{n-1}{l}{m}\enom{n}{l+i}{m}\\
&\stackrel{\mbox{\tiny Id. (ii)}}=&   n \sum_{i=1}^m i \sum_{l=0}^{mn-i} \enom{n-1}{l}{m}\enom{n}{mn-l-i}{m}\\
 &\stackrel{\mbox{\tiny Id. (iv)}}=&  n \sum_{i=1}^m i \enom{2n-1}{mn-i}{m}\\ &\stackrel{\mbox{\tiny Id. (iii)} }=& \frac{mn}{2} \enom{2n}{mn}{m}. \end{eqnarray*}} The third equality is proved in a similar way. \ep
\begin{identity}[Extended Entry 3.81 in~\cite{gould}] \label{3.81}For $n \geq 0$ \[ \arraycolsep=1pt\def\arraystretch{0.5}
\sum_{k=0}^{mn} (-1)^k \enom{n}{k}{m}^2 = \left\{
  \begin{array}{lllll}
    &&0, && \hbox{if $mn$ is odd ;} \\
    && &&\\
    &&\displaystyle \enom{n}{mn/2}{m}, && \hbox{if $m$ is even;} \\
    && &&\\
    &&\displaystyle \sum_{i=0}^n (-1)^i {n \choose i} \enom{2n}{mn/2- i}{\frac{m-1}{2}}, && \hbox{If $m$ is odd and $n$ even.}
  \end{array}
\right.\] \eid \bp Let $A(n)$ denote the LHS of this identity. First, we observe that \[A(n) = \sum_{k=0}^{mn} (-1)^k \enom{n}{k}{m} \enom{n}{mn-k}{m}= [t^{mn}]\left(p_m(-t) \, p_m(t) \right)^n .\] Since the polynomial $\left(p_m(-t) \, p_m(t) \right)^n$ is an even function, $A(n)$ vanishes if $mn$ is odd.

If $m$ is even, the equality $p_m(-t) \, p_m(t) = p_m(t^2)$ does hold. Thus \[ A(n)=[t^{mn}] \sum_{k=0}^{mn} \enom{n}{k}{m} t^{2k}= \enom{n}{mn/2}{m} .\]

If $m$ is odd and $n$ is even, we use the equality $p_m(-t) \, p_m(t) = (1-t^2)\, p_{\frac{m-1}{2}} (t^2)^2$ and write \be \label{last} A(n)= [t^{mn}] \sum_{i,j} (-1)^i \binom{n}{i} \enom{2n}{j}{\frac{m-1}{2}} t^{2i+2j} = \sum_{i=0}^n (-1)^i \binom{n}{i} \enom{2n}{mn/2 - i}{\frac{m-1}{2}} , \ee proving Identity~\ref{3.81}. \ep
 \begin{identity} \label{4nom} For $r \geq 0$, \[\sum_{k=0}^{6r} (-1)^k \enom{2r}{k}{3}^2 = (-1)^{r} {4r \choose r} \, ,\]\eid \bp Write the last term of Eq.~\eqref{last} for $m=3$ and $n=2r$ as \[ \sum_{i=0}^{2r} (-1)^i \binom{2r}{i} \binom{4r}{r +i}. \] Employing the formula (Entry 3.57 in~\cite{gould}) :
\[\sum_{i=0}^{2r} (-1)^i \binom{2r}{i} \binom{2r+2x}{i+x} = (-1)^r \binom{2r}{r} \frac{\binom{2r+2x}{x}}{\binom{r+x}{r}},\] with $x=r$, Identity~\ref{3.81} simplifies to Identity~\ref{4nom}. \ep
\begin{identity}[Extended Identity 155 in~\cite{benj}] \label{idB} For all $n \geq 0$,
 \bea \nonumber \sum_{l =\lceil k/m \rceil }^{n} {n \choose l} \enom{l}{k}{m} &=& \sum_{j=0}^n 2^{n-j} {n \choose j} \enom{j}{k-j}{m-1},  \\
\nonumber \sum_{l =k }^{mn}\enom{n}{l}{m} {l \choose k} &=& \sum_{j=0}^n (-1)^{n-j} {n \choose j} {(m+1)j \choose k+n}. \eea \eid  \bp The GF with respect to $k$ of both sides of the first equation is straightforwardly given by $(1+p_m(t))^n$. Moreover, the GF of the LHS of the second equation is $\left(p_m(1+t)\right)^n$. On the other hand, \begin{eqnarray*}
     [t^k] p_m(1+t)^n =[t^k] t^{-n}\left((1+t)^{m+1}-1\right)^n &=& [t^{k+n}] \sum_{j=0}^n \binom{n}{j}(-1)^{n-j} (t+1)^{j(m+1)}\\
& = & \sum_{j=0}^n (-1)^{n-j} {n \choose j} {(m+1)j \choose k+n},
\end{eqnarray*} proving Identity~\ref{idB}. \ep
\subsection{Some non-trivial identities for trinomial coefficients} The trinomial coefficients ($m=2$) can be written in terms of the Gegenbauer polynomials\footnote{The Gegenbauer polynomials $C_n^{(\alpha)}(x)$ are defined by $(1-2xt+t^2)^{-\alpha}=\sum_{n=0}^\infty C_n^{(\alpha)}(x) t^n$\cite[p.783]{HMF}.} as \be \label{gegen} \enom{n}{k}{2}=C_k^{(-n)}(-1/2)=(-1)^k C_k^{(-n)}(1/2).\ee
Through this connection, many properties of  Gegenbauer polynomials can be specialized to give non-trivial identities for trinomial coefficients. We list some of them with the belief that they may stimulate further extensions.
\begin{identity}[Dilcher formula~\cite{dilcher}] \label{dilcher} If $n \geq 1$, then \[\enom{-n}{k}{2}=(-1)^k \sum_{1 \leq j_1 < \ldots < j_k \leq n+k-1} \left( 1+2 \cos \frac{j_1 \pi}{n+k}\right)\cdot \cdots  \cdot\left( 1+2 \cos \frac{j_k \pi}{n+k}\right).\] \eid
\bp Identity~\ref{dilcher} is a specialization of a formula for Gegenbauer polynomials due to Dilcher~\cite{dilcher}. \ep \begin{remark} The original formula in~\cite{dilcher} contains a minor typographical error. The correct one is \[C_n^{(k)}(x)=2^n \sum_{1 \leq \, j_1 < \ldots < j_n \, \leq n+k-1} \left( x+ \cos \frac{j_1 \pi}{n+k}\right)\cdot \cdots  \cdot\left( x+ \cos \frac{j_k \pi}{n+k}\right), \quad k \geq 0.\] \end{remark}
\begin{identity}[see \cite{rain}, formula 32, p.282] \label{24} For $p > 0$, \[\sum _{k=0}^n \frac{(-1)^k (k+p)}{(n-k)! \, (k+n+1)! \, \displaystyle \binom{2 p+k+n}{2p-1}} \enom{-p}{k}{2} = \frac{1}{2  n! \left(p +1/2\right)_n } \left(\frac{3}{4}\right)^n \, ,\] where the Pochhammer symbol $(x)_n$ is defined by $(x)_n= x(x+1)\cdots (x+n-1)$. \eid
\begin{identity}[see \cite{rain}, formula 36, p.283] For $p > 0$, \[\sum _{k=0}^{\lfloor n/2 \rfloor} \frac{p+n-2k}{k! \, (p)_{n+1-k}}\enom{-p}{n-2k}{2} = \frac{(-1)^n}{n!}\, .\] \eid
 \begin{identity}[A series of products of three trinomial coefficients~\cite{braf}] \label{braf} For $n \geq 1$, \[\sum_{k=0}^\infty (-1)^k \frac{(k+n)\, k!\,^2  }{(k+2n-1)!\,^2} \enom{-n}{k}{2}^3 = \frac{2^n \sqrt{3}}{3^{3n-1}}\frac{\pi}{(n-1)!\,^4} \, .\]   \eid
 \begin{identity}[see \cite{rain}, formula 31, p.282] \label{HGF3} For $n > 0$ and $|t|<1$,  \[\sum_{k=0}^\infty  \frac{(n+1/2)_k }{(2n)_k}\enom{-n}{k}{2}\, t^k=\frac{2^{n-1/2}}{\sqrt{1+t+t^2}}\, \left(1+ \frac{t}{2} + \sqrt{1+t+t^2} \right)^{-n + 1/2}\, .\] \eid
For more formul\ae\ involving Gegenbauer polynomials, we refer the reader to the Dougall's linearization formula~\cite[p.39]{askey}, the Gegenbauer connection formula~\cite[p.59]{askey} and to others in~\cite[Eq. (1.4)]{alam},~\cite[Eq.(19)]{ruiz},~\cite{weisner}. It would be instructive to generalize the above trinomial identities for arbitrary order $m$.

\section{Concluding remarks} \label{s5} In this paper, we have proved a selection of extended binomial coefficient identities; a more ambitious task is obviously to extend all the known properties of the Pascal triangle, as suggested by Comtet~\cite{Comt}.

It is well-known that~\cite[p. 77]{Comt}
\[ \enom{n}{k}{m} = \frac{2}{\pi}\int_0^{\pi/2}
\left(\frac{\sin (m+1) t}{\sin t}\right)^n \cos\left((n m - 2k)t\right) \; dt
\, .\]
Through this integral representation, one possibly may establish combinatorial identities in the manner of Egorychev~\cite{egor}. A natural extension for real values of $n$, $k$ and $m$ can also be considered.

We stress, finally, that several combinatorial interpretations for the polynomial coefficients are known~\cite{fahssi};  an instructive exercice would be to seek, \`a la Benjamin and Quinn~\cite{benj}, combinatorial proofs of the above identities.


\end{document}